\documentclass[11pt]{article}
\usepackage{mathptmx}

\usepackage{epsfig,amsfonts,latexsym,amssymb,shadow,amsmath,wrapfig}
\usepackage{subfigure}
\usepackage{epstopdf}
\usepackage{verbatim}
\usepackage{mathtools}

\usepackage{graphicx,color,multirow,colortbl,hhline,array,url,fancyhdr,setspace,enumerate}
\usepackage[normalem]{ulem}
\usepackage[margin=1.2in]{geometry}
\usepackage{enumitem}

\usepackage{amsmath,bm}
\usepackage{mathtools} 
\usepackage{float}
\usepackage[ruled,vlined,linesnumbered]{algorithm2e}
\usepackage{booktabs}
\usepackage{tabularx}

\usepackage{threeparttable}

\usepackage{bm}

\usepackage{xcolor}

\usepackage[super,sort&compress,comma]{natbib} 

\usepackage[hang,flushmargin,symbol]{footmisc}

\usepackage{bigstrut}
\usepackage{color, colortbl}
\makeatletter
\newlength\mylena
\newlength\mylenb
\newcommand\mystrut[1][2]{%
    \setlength\mylena{#1\ht\@arstrutbox}%
    \setlength\mylenb{#1\dp\@arstrutbox}%
    \rule[\mylenb]{0pt}{\mylena}}
\makeatother
\newcolumntype{L}{>{\centering\arraybackslash}m{3cm}}
\newcolumntype{K}{>{\centering\arraybackslash}m{1.2cm}}
\newcolumntype{J}{>{\centering\arraybackslash}m{1.6cm}}
\newcolumntype{M}{>{\centering\arraybackslash}m{3cm}}
\newcolumntype{N}{>{\centering\arraybackslash}m{4cm}}
\newcolumntype{O}{>{\centering\arraybackslash}m{3cm}}

\begin{document}
	
\title{Sector coupling via hydrogen to lower the cost of energy system decarbonization}
\begingroup

\author{Guannan He\footnotemark[1]~~\footnotemark[2]~~,
    Dharik S. Mallapragada\footnotemark[2]~~,
	Abhishek Bose\footnotemark[2]~~,\\
	Clara F. Heuberger\footnotemark[3]~~, 
	Emre Gen\c{c}er\footnotemark[2]
}

\footnotetext[1]{Corresponding author: gnhe@mit.edu}
\footnotetext[2]{MIT Energy Initiative, Massachusetts Institute of Technology, Cambridge, MA, USA.}
\footnotetext[3]{Shell Global Solutions International B.V., Shell Technology Centre Amsterdam, 1031 HW Amsterdam, Netherlands.}

\endgroup

\date{}

\maketitle

\noindent \textbf{Abstract}:  
There is growing interest in hydrogen (H$_\text{2}$) use for long-duration energy storage in a future electric grid dominated by variable renewable energy (VRE) resources. Modelling the role of H$_\text{2}$ as grid-scale energy storage, often referred as “power-to-gas-to-power (P2G2P)” overlooks the cost-sharing and emission benefits from using the deployed H$_\text{2}$ production and storage assets to also supply H$_\text{2}$ for decarbonizing other end-use sectors where direct electrification may be challenged. Here, we develop a generalized modelling framework for co-optimizing energy infrastructure investment and operation across power and transportation sectors and the supply chains of electricity and H$_\text{2}$, while accounting for spatio-temporal variations in energy demand and supply. Applying this sector-coupling framework to the U.S. Northeast under a range of technology cost and carbon price scenarios, we find a greater value of power-to-H$_\text{2}$ (P2G) versus P2G2P routes. P2G provides flexible demand response, while the extra cost and efficiency penalties of P2G2P routes make the solution less attractive for grid balancing. The effects of sector-coupling are significant, boosting renewable energy generation by 12-55\% with both increased capacities and reduced curtailments and reducing the total system cost (or levelized costs of energy) by 6-14\% under 96\% decarbonization scenarios. Both the cost savings and emission reductions from sector coupling increase with H$_\text{2}$ demand for other end-uses, more than doubling for a 96\% decarbonization scenario as H$_\text{2}$ demand quadraples. Moreover, we found that the deployment of carbon capture and storage is more cost-effective in the H$_\text{2}$ sector because of the lower cost and higher utilization rate. These findings highlight the importance of using an integrated multi-sector energy system framework with multiple energy vectors in planning energy system decarbonization pathways.

\section{Introduction}
As the greenhouse gas (GHG) emission intensity of electricity generation in various regions has declined with continued adoption of wind and solar generation, there is growing interest to pursue electrification-centric decarbonization strategies for other end-use sectors where emissions reduction has been sluggish. Yet, direct electrification may be practically challenged for some of these end-uses, such as in the case of heavy-duty transport where volumetric energy density and refueling time are key drivers for fuel choice. In this context, there is renewed interest in hydrogen (H$_\text{2}$) and H$_\text{2}$ derived carriers for their role in the decarbonization of difficult-to-electrify end-uses in transport, building and industrial sectors. In addition to the plurality of its end-uses, the multiple technology choices across the H$_\text{2}$ supply chain, from production, storage, transport and end-use, make its assessment a complex systems problem. 
Here, we propose a scalable decision-support framework for assessing the impact of technology and policy choices on the decarbonization of power sector in conjunction with other end-use sectors. This framework provides a systematic way to study the role and impact of H$_\text{2}$-based technology pathways in a future low-carbon, integrated energy system at a regional/national scale.

Recent renewed interest in H$_\text{2}$ has been partially intrigued by expectations of a future renewables-dominant electric grid and cost declines for water electrolyzers\cite{Schmidt17}, both of which raise the prospect of electrolytic H$_\text{2}$ becoming cost-competitive with fossil fuel-based pathways, such as natural gas reforming\cite{GUERRA20192425,MALLAPRAGADA2020100174}. Besides the economics of electrolytic H$_\text{2}$ production \cite{Finke2021,Glenk19}, many studies have focused on evaluating the economics of H$_\text{2}$-based energy storage (power-to-gas-to-power, P2G2P), which relies on electrolysis for H$_\text{2}$ production, under deep decarbonization scenarios. Some of the studies in this area focus on: 
1) comparing the cost-effectiveness of P2G2P with other types of long-duration energy storage options like pumped hydro and compressed air energy storage for VRE integration, from a marginal deployment perspective (i.e. electricity price taker) \cite{Guerra2020,Clerjon19,Pellow15}, 2) assessing least-cost investment and operation of H$_\text{2}$ storage and short-term energy storage like lithium-ion batteries in the context of a VRE dominant power systems \cite{Pellow15,Dowling2020,ZERRAHN20171518,MCPHERSON2018649}, and 3) the operational scheduling of H$_\text{2}$ storage in power markets \cite{el-taweel19,yang19}. Although these studies provide useful insights to compare different energy storage technologies from the perspective of the power sector, they overlook the multiple potential uses of H$_\text{2}$ (or H$_\text{2}$ derived carriers) outside the power sector and the associated cost-savings resulting from sharing infrastructure costs across these uses. Consequently, in the absence of modelling sector-coupling interactions, the role of H$_\text{2}$ storage may be under-valued as compared to other long-duration storage technologies in future low-carbon power grids\cite{David_science_2018}.

With the above motivation, a number of studies have expanded the scope of traditional power sector capacity expansion models (CEM) to endogenize investment decisions in end-use technologies, which includes some parts of the H$_\text{2}$ supply chain, notably electrolytic H$_\text{2}$ production. These studies highlight the potential for flexible electricity consumption in other end-uses to partially substitute the need for energy storage in the electricity sector and alter generation mix in the power sector towards increasing VRE deployment \cite{D0EE02016H,VICTORIA2019111977,BODAL202032899,BROWN2018720}. While these studies are inspiring, the interactions between the H$_\text{2}$ supply chain and the power sector, in many of the studies, exclude critical components in the H$_\text{2}$ supply chain. For example, some studies ignore the possibility of natural gas-based H$_\text{2}$ production from steam methane reformer (SMR) with or without carbon capture and storage (CCS) \cite{li_19,reu17,welder_18}, which is the dominant mode of H$_\text{2}$ production today. Second, most literature either do not consider some modes of H$_\text{2}$ transmission\cite{li_19,ochoa18,de-leon_almaraz_hydrogen_2014} or when it is included, the modelling of H$_\text{2}$ transmission is oversimplified by setting fixed lower and upper H$_\text{2}$ flow limits for each route \cite{reu17,welder_18,ochoa18,li20}. These approaches may not capture the potential benefits of both H$_\text{2}$ pipeline and trucks serving as transmission and storage assets simultaneously. Notably, H$_\text{2}$ trucks can function as mobile storage, which has been shown to provide greater operational flexibility than stationary storage\cite{he2020hsc}. Moreover, the existing literature does not reveal a clear evolution of the role of H$_\text{2}$ in energy systems as the costs of H$_\text{2}$ infrastructure decline with increased adoption or technology innovation. 

This paper develops a high-fidelity electricity-H$_\text{2}$ (e-H$_\text{2}$) capacity planning model to study the role of H$_\text{2}$ in low-carbon energy systems, the sector-coupling effects, and the trade-offs between various technology options across the entire bulk supply chain\footnote{The last-mile distribution networks of electricity and H$_\text{2}$ are not considered in the e-H$_\text{2}$ model.} of both energy carriers. For a pre-defined set of electricity and H$_\text{2}$ demand scenarios, the model determines the least-cost technology mix across the power and H$_\text{2}$ sectors while adhering to operational constraints of the power and H$_\text{2}$ supply chains at an hourly resolution along with the spatiotemporal variations in VRE supply and energy demands. 

Applying the e-H$_\text{2}$ model to the U.S. Northeast energy system for a range of CO$_\text{2}$ price (up to \$1000/tonne CO$_\text{2}$), H$_\text{2}$ demand and technology cost scenarios, we find electrolytic H$_\text{2}$ supply to be cost-effective under moderate carbon policy (\$50/tonne or greater) and/or electrolyzer capital cost of \$500/kW or lower. The interactions between the power and H$_\text{2}$ supply chains increase investments in VRE generation and reduce investments in dispatchable resources like battery storage and natural gas generation in the power sector, which results in reducing the total system cost by up to 14\% in the most carbon constrained scenarios analyzed here. Notably, as opposed to most literature that emphasizes the role for H$_\text{2}$ in the power sector as grid-scale energy storage, i.e. P2G2P\cite{Guerra2020,Guerra2020,Pellow15,Clerjon19,Dowling2020,el-taweel19,yang19}, we find a greater role for H$_\text{2}$ to serve as a flexible demand response resource based on the use of electrolysis in conjunction with H$_\text{2}$ storage. This finding, stemming from the reduced capital cost and energy efficiency of the P2G path, is found to be 
valid across a range of CO$_\text{2}$ price scenarios as well as capital cost assumptions for electrolysis and G2P generation. Finally, we find that the role for natural gas in a future energy system is predominantly in the H$_\text{2}$ supply chain and this role remains robust to increasing CO$_\text{2}$ prices because of the relative cost-competitiveness of CCS-equipped natural gas based H$_\text{2}$ production vs. electrolytic H$_\text{2}$ production.

	
\section{Methods}
\subsection{Model overview}

\begin{figure}[!htb]
	\centering
	\includegraphics[width=\textwidth]{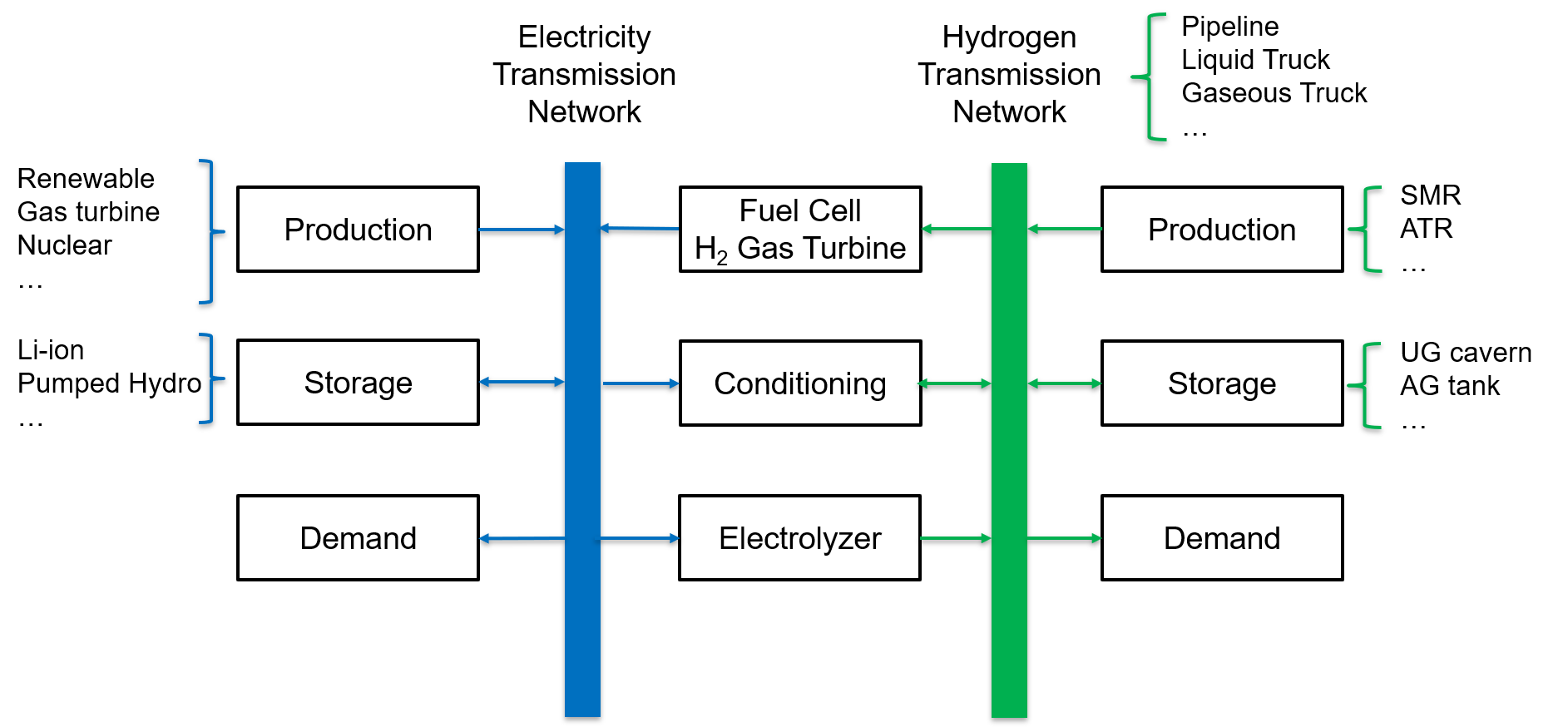}
	\caption{Superstructure of the coupled model of power systems and H$_\text{2}$ supply chain.}
	\label{fig-schematic}
\end{figure}

The developed e-H$_\text{2}$ planning model evaluates investments and operations across the bulk supply chain for electricity and H$_\text{2}$, including production, storage, transmission, conditioning (compression/liquefaction in the case of H$_\text{2}$) and demand as shown in Figure \ref{fig-schematic}. The model determines the least-cost mix of electricity and H$_\text{2}$ production, storage, and transmission infrastructures to meet power and H$_\text{2}$ demands subject to a variety of operational and policy constraints. The developed model can incorporate a wide range of power and H$_\text{2}$ technology options, including VRE generation, carbon capture and storage (CCS) applied to power and H$_\text{2}$ generation, and truck (gaseous and liquid) and pipelines for H$_\text{2}$ transportation. The power systems and H$_\text{2}$ supply chain are coupled through electrolysis and power generation technologies fueled by H$_\text{2}$, as well as electricity consumption in H$_\text{2}$ compression/liquefaction. 
The operational constraints of the model, implemented at an hourly resolution, include: a) supply-demand balance for H$_\text{2}$ and electricity at each zone, b) inventory balance constraints for stationary storage technologies, c) inventory balance constraints related to trucks at a given location (any of the zones and routes, arriving, departing or in transit) and for different states (empty and full), and d) linearized unit commitment for conventional thermal power generation technologies and natural gas based H$_\text{2}$ production technologies. We model these operational constraints at an hourly resolution over a set of representative weeks that are selected from applying time-series clustering to annual demand and VRE resource profile data\cite{MALLAPRAGADA2020115390}, to approximate annual system operations. The time-domain reduction preserves chronological variability of  power and H$_\text{2}$ demands and VRE resource availability and the correlations among them, while reducing the model size to still be computationally tractable. Process level CO$_\text{2}$ emissions are penalized with a price on emissions that is applied to both sectors. The details of the power system planning model and the H$_\text{2}$ supply chain model can be referred to \citet{genx} and \citet{he2020hsc}, respectively. The codes and data are available on GitHub \cite{gn_he_2021_4443233}. 

\subsection{Case study setup}
\begin{figure}[!htb]
    \centering
    \includegraphics[width =
    \textwidth/2]
    {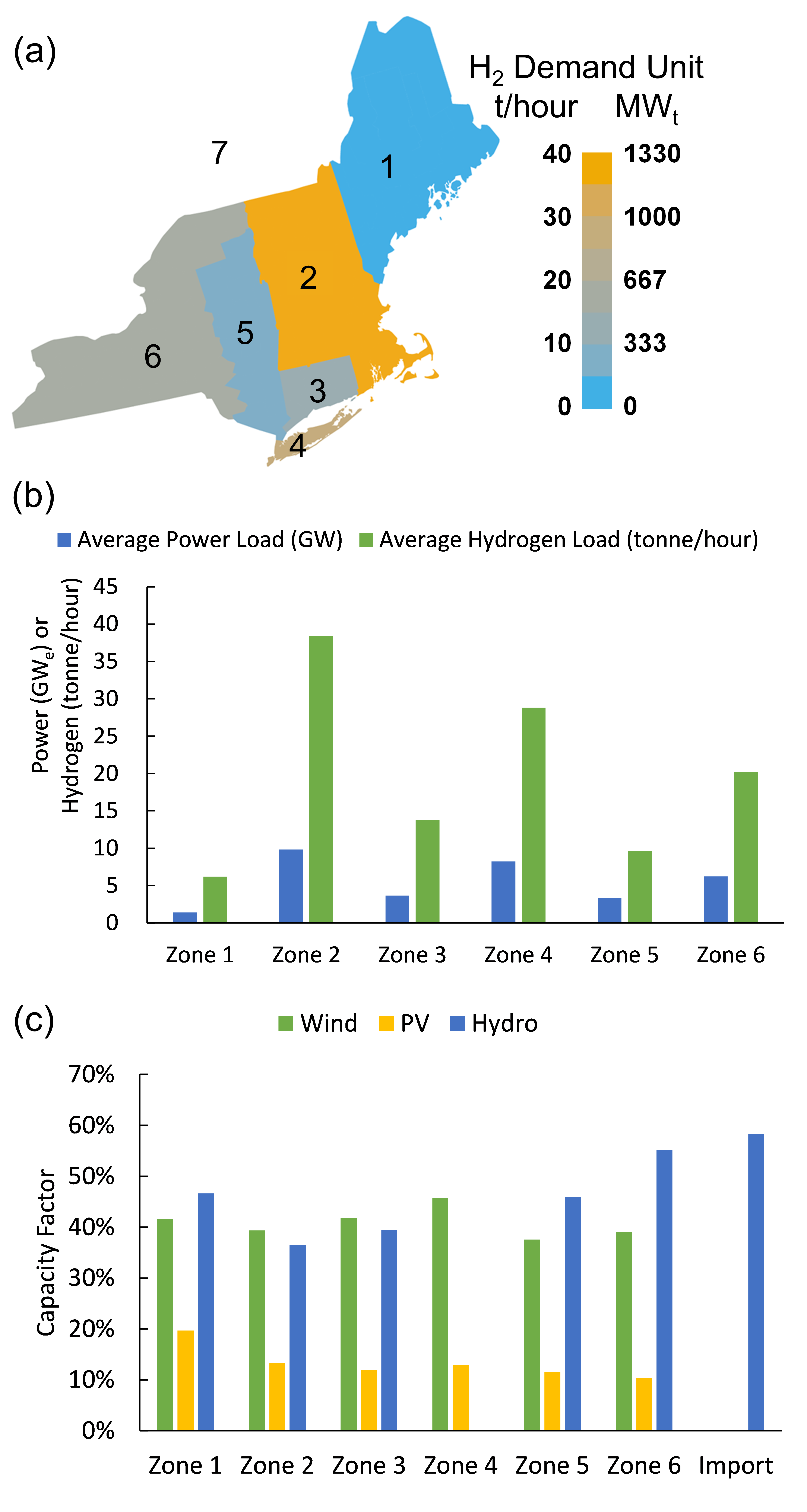}
    \caption{Demand and renewable energy resources distributions in the U.S Northeast. H$_\text{2}$ demands are estimated based on transportation fuel consumption data in 2017, and electricity demands are projections for 2050. (a) Geographical zone classification for U.S. North-East and average H$_\text{2}$ demands for each zone; (b) Average power and H$_\text{2}$ demands for each zone; (c) Average capacity factors of wind, PV, and hydro for each zone.}
    \label{fig:region_setting}
\end{figure}

We illustrate the value of the proposed model using a case study where we assess electricity and H$_\text{2}$ infrastructure outcomes for the U.S. Northeast region under a variety of demand, technology and CO$_\text{2}$ price scenarios for 2050. We model a greenfield 2050 system with the exception of existing inter-zonal transmission and hydro power generation (both domestic and imports from Canada) and pumped hydro storage capacity in the region. The U.S. Northeast region is represented in the model as six zones, shown in Figure \ref{fig:region_setting}(a), according to the zonal boundaries adopted from the Integrated Planning Model \cite{epa2020v}. An additional seventh zone is included with zero energy demand to represent imports of Canadian hydro power generation that is limited by power and transmission capacity constraints. States in the Independent System Operator New England (ISO-NE) are split into zone 1-3 and the New York Independent System Operator (NYISO) is split into zone 4-6 based on their load share split in 2012. As zone 4 is heavily urbanized, we do not allow centralized H$_\text{2}$ generation (SMR or SMR with CCS) to be built in that zone, while distributed electrolyzers are allowed.


Electricity demand data (excluding electrolysis) are based on 2018 NREL electrification futures study load projection for 2050 \cite{mai2018electrification}, with assumed business-as-usual technology advancement and reference electrification. 
The H$_\text{2}$ demands for each zone are developed based on available fuel consumption data and hourly refueling profiles \cite{nexsnt} for both light- and heavy-duty fuel cell electric vehicles (FCEV) and the relative penetration of FCEV. In the base case, we assume a 20\% FCEV penetration for the transport sector in the U.S Northeast. Light-duty vehicle (LDV) fuel consumption for each zone is estimated using state-level gasoline consumption data for 2017 from the the U.S. Energy Information Administration, which is then converted to a H$_\text{2}$ consumption equivalent based on the relative efficiency of FCEV to gasoline internal combustion vehicle. Heavy-duty vehicle (HDV) demand projections are based on the National Freight Analysis Framework \cite{faf}.
The zonal average demands of power and H$_\text{2}$ are shown in Figure \ref{fig:region_setting}(b).

For the power system, we include thermal, renewable, nuclear generation, and storage resources, whose main parameters are derived from the NREL annual technology baselines\cite{vimmerstedt20182018} and the EIA Annual
Energy Outlook 2018\cite{us2018annual} for the year 2050, as summarized in Table \ref{tab:Power_Gen}. As no new coal plants will likely be built in this region, we consider the options of natural gas combined cycle gas turbine (CCGT), open cycle gas turbine, as well as CCGT with CCS. The VRE resource cost and availability in each zone are represented by supply curves \cite{BROWN2020} to characterize different possible sites with specific resource profile, maximum potential capacity, and average cost of interconnection. We use three supply curves per zone for onshore wind and one supply curve per zone for PV. Offshore wind is included with no capacity limits and single resource profile for zone 2 and zone 4 based on sampling sites from the NREL Wind Toolkit that overlaps with the areas. Distributed PV is modelled with a separate resource profile per zone and minimum build requirement to meet 2029 projections by NYISO and ISONE. For hydropower, we consider hydro reservoir, hydro run-of-river, and Canadian hydro and extrapolate hourly generation from historical monthly output. The average capacity factors of VRE resources are shown in Figure \ref{fig:region_setting}(c). Lithium-ion battery storage and pumped hydro storage are considered for electrical energy storage. The initial power transfer capacities between each zone are developed from the integrated planning model (IPM) model documentation \cite{epa2020v} and are tabulated in Table S2 in the SI, while transmission expansion costs are listed in Table S1.

\begin{table}[tb]
	\centering
	\caption{Major Parameters for generation and storage technologies in the power sector for the year 2050. CAPEX: capital cost; FOM: fixed operational and maintenance cost; VOM: fixed operational and maintenance cost.}
	\resizebox{\textwidth}{!}{
		\begin{tabular}{ l K K K J K K K K K K K }
			\hline
			
			Technology & Onshore Wind & Offshore Wind & Utility PV & Distributed PV & Li-ion Battery & Pumped Hydro & CCGT & OCGT & CCGT w/ CCS & Nuclear\\
			\hline
			Power CAPEX ($10^3$\$/MW) & 1,086 & 1,902 & 724 & 1,083 & 120 & 1,966 & 817 & 816 & 1798 & 6126\\
			Energy CAPEX ($10^3$\$/MWh) & - & - & - & - & 126 & - & - & - & - & -\\

			FOM ($10^3$\$/MW-year) & 35 & 48 & 11 & 9 & 2 & 41 & 11 & 12 & 34 & 105\\
			
			%
			%
			
			VOM (\$/MWh) & - & - & - & - & - & - & 3 & 7 & 7& 2\\

			Heat Rate (MMBTU/MWh) & - & - & - & - & - & - & 6 & 9 & 7 & 10\\

			Round-trip Efficiency & - & - & - & - & 85\% & 80\% & - & - & - & -\\

			Lifetime (years) & 30 & 30 & 30 & 30 & 15 & 50 & 30 & 30 & 30 & 30\\
			\hline
		\end{tabular}
	}
	\label{tab:Power_Gen}
\end{table}%

\begin{table}[ !h]
	\small{
		\centering
		\caption{Major parameters for H$_\text{2}$ generation and G2P technologies. CAPEX: capital cost.}
		\resizebox{\textwidth}{!}{
			\begin{tabular}{ l c c c c c }
				\hline
				
				& Electrolysis\cite{IEA2019} & SMR\cite{IEA2019} & SMR w/ CCS\cite{IEA2019} & Fuel Cell\cite{Battelle2016} & CCGT-H$_\text{2}$\cite{vimmerstedt20192019}\\
				%
				\hline
				
				Unit CAPEX & 450 \$/kW$_\text{e}$ & 910 \$/kW$_{\text{H}_2}$ & 1,280 \$/kW$_{\text{H}_2}$ & 1,264 \$/kW$_\text{e}$ & 1,171 \$/kW$_\text{e}$\\

				Lifetime (years) & 10 & 25 & 25 & 10 & 25\\
				%

				Efficiency (LHV) & 74\% & 76\% & 69\% & 60\% & 65\%\\
				
				%
				%
				
				Emissions Intensity (tonne $\mathrm{CO_{2}}$/tonne $\mathrm{H_{2}}$) & 0 & 8.9 & 1.0 & 0 & 0\\
				\hline
			\end{tabular}
			\label{tab:H2_Gen}
		}
	}
\end{table}%

\begin{table}[ !h]
	\centering
	\caption{Major parameters for H$_\text{2}$ transmission and storage technologies. CAPEX: capital cost; OPEX: operational cost; A: cost and electricity consumption proportional to pipeline length; B: cost and electricity consumption irrelevant to pipeline length; C: truck and tank storage compression related costs and electricity consumption.}
	\resizebox{\textwidth}{!}
	{
		\begin{tabular}{l c c c c}
			\hline
			
			& Pipeline & Gas Tank & Liquid Truck & Gas Truck\\
			\hline
			
			Unit Capacity & 0.3 tonne/mile & 0.3 tonne & 4 tonne & 0.3 tonne\\

			Capital Cost & 2.8 M\$/mile \cite{FEKETE201510547,ingaa} & 0.58 M\$/unit\cite{yang_determining_2007} & 0.8 M\$/unit\cite{yang_determining_2007} & 0.3 M\$/unit\cite{yang_determining_2007}\\

			Compression CAPEX (A) (\$/mile-unit) & 700\cite{schoenung11,samsatli_multi-objective_2018} & 0 & 0 & 0\\

			Compression CAPEX (B) (\$/unit) & 0.75 & 0 & 0 & 0\\

			Compression Electricity (A) (MWh/tonne-mile) & 1 & 0 & 0 & 0\\

			Compression Electricity (B) (MWh/tonne) & 1 & 0 & 0 & 0\\

			Unit OPEX (\$/mile) & 0 & 0 & 1.5 & 1.5\\

			Compression CAPEX (C) (\$/(tonne/hr))& 0 & 0.5\cite{yang_determining_2007} & 32\cite{yang_determining_2007} & 1.5 \cite{yang_determining_2007}\\

			Compression Electricity (C) (MWh/tonne) & 0 & 2 \cite{samsatli_multi-objective_2018,schoenung11} & 11 \cite{samsatli_multi-objective_2018,schoenung11} & 1 \cite{samsatli_multi-objective_2018,schoenung11}\\

			Boiloff Rate & 0 & 0 & 3\% & 0\\
			\hline
			
		\end{tabular}
	}
	\label{tab:hsc}
\end{table}%

The main cost and performance parameters of H$_\text{2}$ generation and G2P technologies are summarized in Table \ref{tab:H2_Gen}, which include electrolysis and natural gas fueled SMR, with and without CCS (90\% capture), stationary fuel cell, and H$_\text{2}$ fueled CCGT. We model trucks and pipelines as the key modes of H$_\text{2}$ transmission, with the distance traveled in each case measured by the distances between the polygon centroids of each zone. At the same time, we also model them as potential storage resources, in tandem with stationary H$_\text{2}$ storage. We model the potential deployment of two types of trucks, based on handling H$_\text{2}$ as a cryogenic liquid or compressed gas, while the pipelines are considered as multiples of an 8" pipeline being built across different geographies. We do not consider geological H$_\text{2}$ storage as its availability in U.S Northeast region is uncertain \cite{ug_stor}. The parameters of transmission and storage technologies are summarized in Table \ref{tab:hsc}. The interfaces of each of these transmission and storage technologies with H$_\text{2}$ generation and demand require compression and/or liquefaction depending on the state of H$_\text{2}$. The compression/liquefaction costs comprise of the capital cost of the equipment as well as the operational costs from electricity consumption, which are also provided in Table \ref{tab:hsc}.

\subsection{Modelling approximations}
To maintain computational tractability with the expanded scope of investment and operational decisions considered here, we implement the following approximations, whose potential impacts on model outcomes are described below. Similar to other power sector CEM studies \cite{Guerra2020,BROWN2020}, we approximate annual hourly system operations based on modelling operations of the system over 30 representative weeks, selected using K-means based clustering techniques, described elsewhere \cite{MALLAPRAGADA2020115390}, in conjunction with heuristics regarding so-called "extreme" weeks applied to 7-years (2007-2013) of load and VRE availability data. Because representative periods identified via clustering techniques are known to emphasize typical weeks over extreme weeks, we apriori identified weeks in the data set with the highest load and lowest average capacity factors and added them to the set of 30 representative weeks to be considered by the model.

Some of the technologies considered in the power and H$_\text{2}$ supply chain, namely thermal power plants, H$_\text{2}$ pipelines and SMR based H$_\text{2}$ production facilities, exhibit economies of scale and limited operational flexibility, which typically requires using binary or integer variables to represent their investment and operations. However, because in nearly all scenarios, we are deploying more than one unit of each technology, the approximation of modelling investment in these technologies as continuous rather integer variables is relatively small. Prior modelling work \cite{palmintier2013incorporating,6684593} has shown that such an approximation in practice results in a relatively small error in the overall dispatch and objective function while leading to large reductions in computational run times. In the case of SMRs, which represent centralized H$_\text{2}$ production sources, 
we estimated that the additional H$_\text{2}$ storage cost (pressurized gas tank) needed to make SMR output flexible only accounts for approximately 2\% of the capital cost of SMR with CCS\footnote{Assuming that a 0.5 hour gas storage is installed to buffer SMR starting up or shutting down, then the capital cost of storage for SMR of one tonne/hour is \$1.9 million (\$0.58 million for a 0.3 tonne unit), which is 2\% of the capital cost of SMR with CCS (\$391 million for a 9.2 tonne/hour unit). The electricity operating cost is negligible when SMR is ramping down (charging and compressing), as the electricity must be very cheap at the time.}
. Therefore, we assume SMR is as flexible as electrolysis in this study.





\section{Results}

\subsection{Optimal technology mixes}
\begin{figure}[ !h]
	\centering
	\includegraphics[width=\textwidth/7*6]{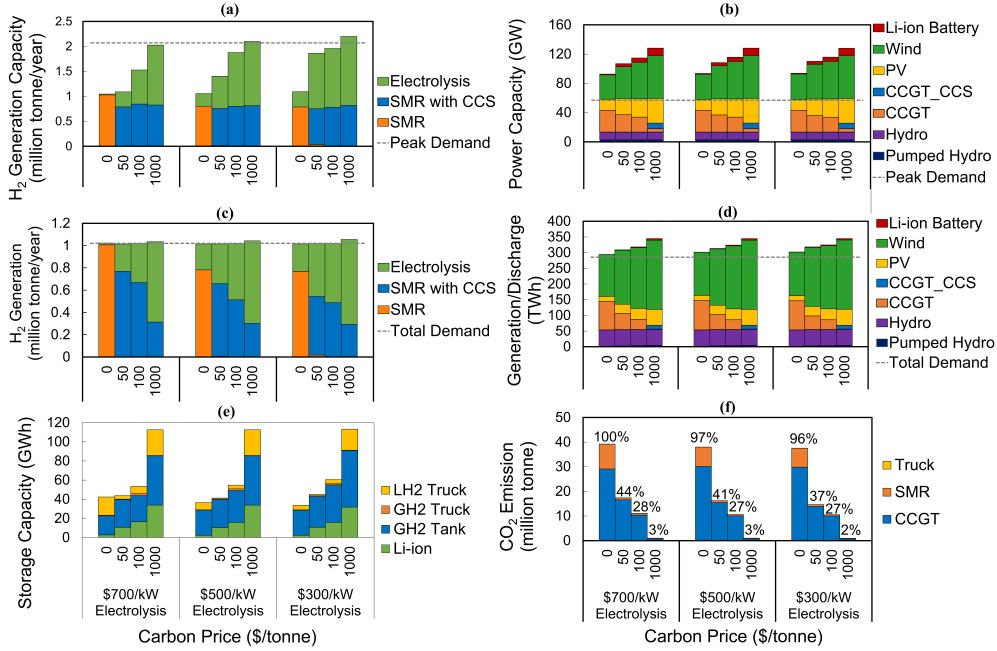}
	\caption{Optimal generation and storage capacity mixes and CO$_\text{2}$ emissions in the power and H$_\text{2}$ sectors under various CO$_\text{2}$ price and electrolysis cost scenarios in U.S. Northeast. Technologies that are not cost-competitive are not shown. (a) H$_\text{2}$ generation capacity; (b) Power generation capacity; (c) H$_\text{2}$ generation per year; (d) Electricity generation or electrical storage energy discharge per year; (e) Energy storage capacity in power and H$_\text{2}$ sectors; (f) CO$_\text{2}$ emissions in power and H$_\text{2}$ sectors. SMR: steam methane reformer; LH2: Liquid H$_\text{2}$; GH2: Gaseous H$_\text{2}$.}
	\label{fig-gen_cap_mix}
\end{figure}

Fig. \ref{fig-gen_cap_mix} presents the optimal technology mixes with different CO$_\text{2}$ prices and electrolyzer capital costs. The scenarios with 0, \$50/tonne, \$100/tonne, and \$1000/tonne CO$_\text{2}$ prices correspond to approximately 100\%, 40\%, 30\%, and 3\% CO$_\text{2}$ emissions compared to no CO$_\text{2}$ price scenario and represent no carbon policy, moderate carbon policy, and deep decarbonization scenarios, respectively. 
Here we highlight several observations from Fig. \ref{fig-gen_cap_mix}:  First, from Fig. \ref{fig-gen_cap_mix} (a) and (b), we can find that as the CO$_\text{2}$ price increases, the H$_\text{2}$ generation switches from central SMR, to SMR with CCS, and then to electrolyzer, accompanied by the power generation shifting away from CCGT to wind and solar. Second, electrolyzer is only cost-effective for deployment at lower CO$_\text{2}$ prices and/or reduced capital costs compared to 2020 costs levels for multi-MW systems, which is near \$800 to \$1000/kW \cite{IEA2019}. 
Thirdly, although CO$_\text{2}$ price increase favors an increasing share of H$_\text{2}$ generation from electrolyzer and renewable power generation, it has a relatively small impact on the installed capacities of natural-gas-fueled H$_\text{2}$ and power generation (SMR and CCGT). CCGT and SMR with CCS remain cost-effective sources of flexible power and H$_\text{2}$ supply for time periods when the lack of VRE generation result in scarcity pricing in the power system for all CO$_\text{2}$ price scenarios analyzed here (see Fig. \ref{fig-gen_cap_mix} (c) and (d)). 
Fourthly, we see CCS utilized at lower CO$_\text{2}$ prices in the H$_\text{2}$ sector (less than \$50/tonne) than in the power sector (greater than \$100/tonne), comparing Fig. \ref{fig-gen_cap_mix} (c) and (d), which implies that if CCS resource availiability is limited, equipping SMR with CCS will be of higher priority. This finding is a result of lower cost of CO$_\text{2}$ capture at SMR facilities than CCGT power plants \footnote{The levelized CO$_\text{2}$ abatement cost of CO$_\text{2}$ capture facility is approximately \$40/tonne at SMR and \$110/tonne at CCGT, given the cost assumptions in the study.} as well as the higher utilization factor of SMR-CCS facilities vs. CCGT-CCS facilities in the analyzed scenarios. For example, in the scenario with \$300/kW electrolysis and \$100/tonne CO$_\text{2}$ price, the capacity factor is 18\% for CCGT, while 64\% for SMR with CCS.
Lastly, we find from Fig. \ref{fig-gen_cap_mix} (e) that H$_\text{2}$ storage, both stationary and mobile, accounts for the majority of storage resources in no and moderate carbon policy scenarios, while the requirement for electrical storage increases with higher VRE penetration in the power sector. Overall, in the future energy system in the U.S. Northeast, we find that natural gas could play a key role as a flexible resource and electrolytic H$_\text{2}$ supply will be cost-effective with moderate carbon policy (\$50/tonne or greater) and/or electrolyzer capital cost reduction (\$500/kW or lower). 

\subsection{Sector coupling effects}
\begin{figure}[!htb]
	\centering
	\includegraphics[width=\textwidth]{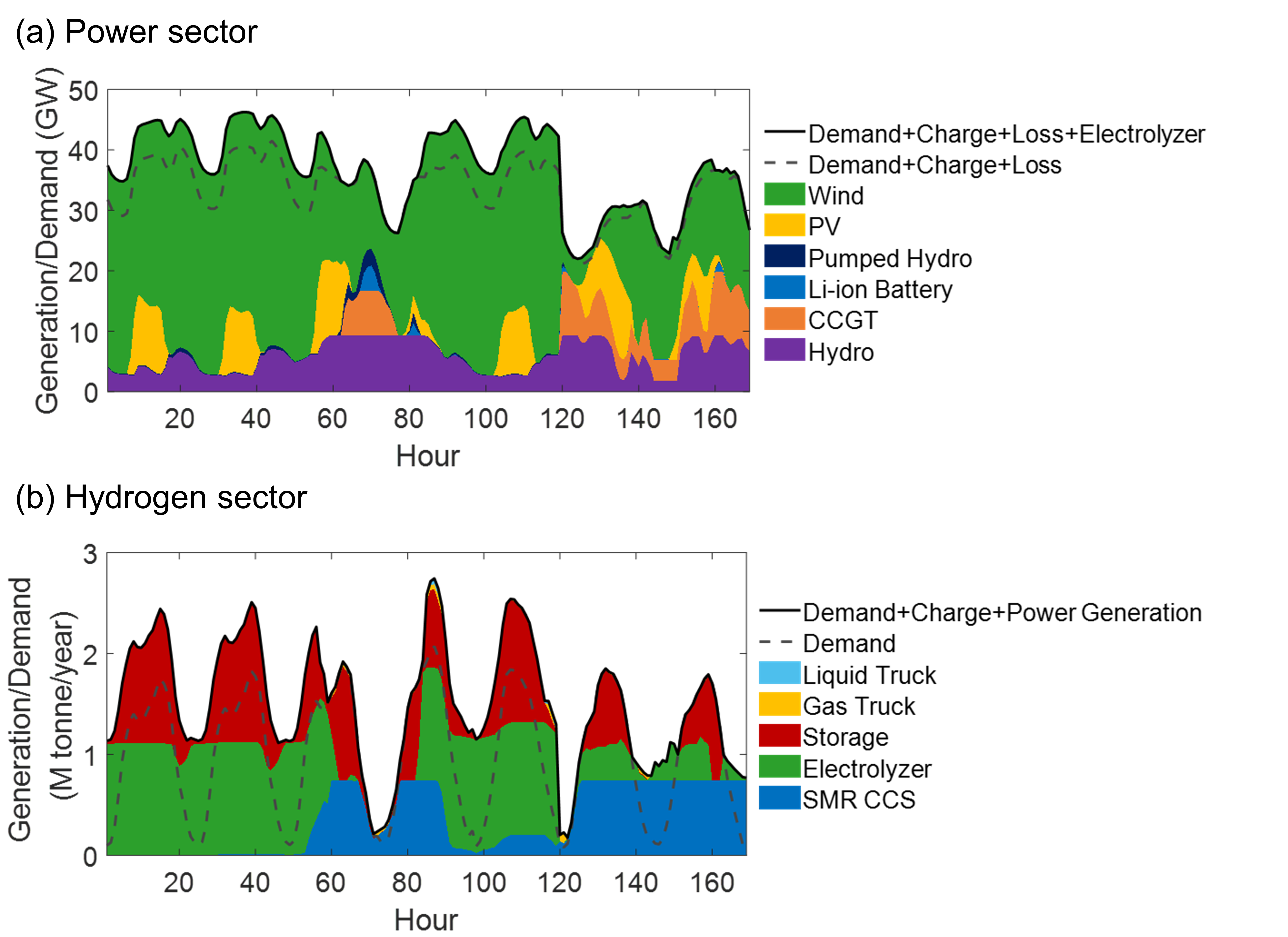}
	\caption{Generation and demand profiles in the power and H$_\text{2}$ sectors in a representative week in the base case demand scenario (H$_\text{2}$ demand corresponding to 20\% FCEV penetration), \$100/tonne CO$_\text{2}$ price and \$300/kW$_\text{e}$ electrolyzer capital cost. (a) Generation and demand profiles in the power sector; (b) Generation and demand profiles in the H$_\text{2}$ sector.}
	\label{fig-gd_profile}
\end{figure}

When the power and H$_\text{2}$ sectors are tightly coupled through electrolysis or H$_\text{2}$-based power generation, the operational flexibility resources in the H$_\text{2}$ sector can support VRE integration in the power sector, leading to overall system cost reductions. Fig. \ref{fig-gd_profile} demonstrates how the two sectors coordinate with each other in a representative week for the scenario with \$300/kW electrolyzer and \$100/tonne CO$_\text{2}$ price. As shown in Fig. \ref{fig-gd_profile}, electrolyzer is the main H$_\text{2}$ supply source when VRE supply is abundant, such as hour 0 to 50. Once VRE is in short supply relative to baseline electricity demand (hour 60 to 80), the SMR, stationary gas storage, and gas and liquid trucks are utilized to meet H$_\text{2}$ demands.

\begin{figure}[!htb]
	\centering
	\includegraphics[width=\textwidth]{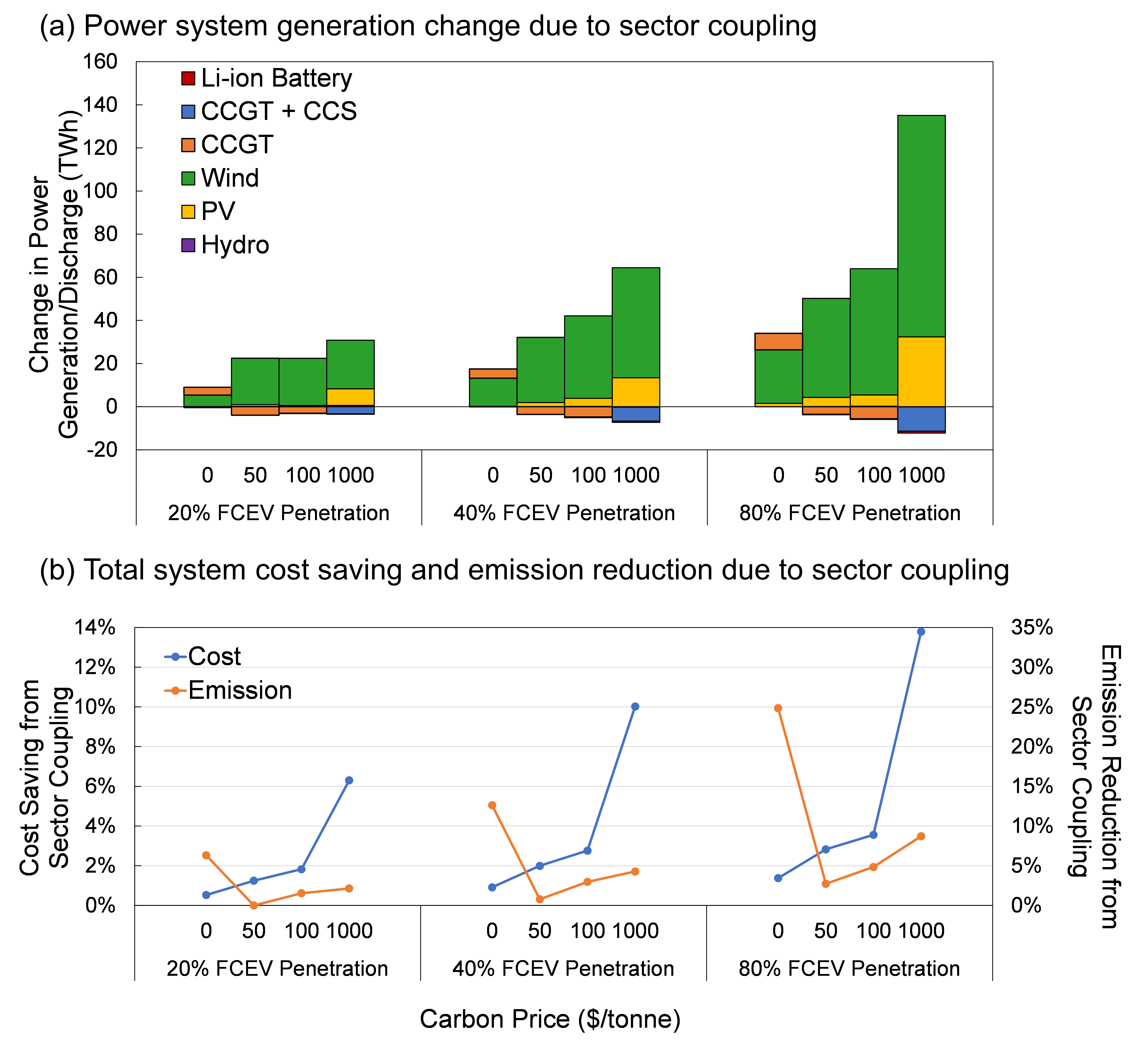}
	\caption{Differences in the optimal power system generation mix, total system cost, and CO$_\text{2}$ emission between energy systems with and without conversions between power to H$_\text{2}$ under various CO$_\text{2}$ price and FCEV penetration scenarios. (a) Differences in the optimal power system generation mix; (b) Differences in the total system cost and CO$_\text{2}$ emission. The capital cost of electrolyzer is assumed to be \$300/kW$_\text{e}$. In (b), the cost savings are shown as percentages of the total system costs in each CO$_\text{2}$ price and FCEV penetration scenario, while the emission reductions are shown as percentages of the CO$_\text{2}$ emission in the case with \$0/tonne CO$_\text{2}$, \$300/kW electrolyzer, and 20\% FCEV penetration.}
	\label{fig-synergy}
\end{figure}

How much is the benefit of coupling power and H$_\text{2}$ sectors? To quantify the impact of sector coupling, we compare the optimal power sector generation mixes with and without the options of conversion between power and H$_\text{2}$ (electrolysis and H$_\text{2}$-based power generation). In the latter case, H$_\text{2}$ production would exclusively rely on natural gas based pathways\footnote{Electricity consumption by conditioning (compression and liquefaction) is still supplied by the power sector, and its cost is accounted in the same way as the coupled case with conversion between power and H$_\text{2}$.}. From Fig. \ref{fig-synergy} (a), we can observe that the power and H$_\text{2}$ interactions (mainly through electrolysis) boost VRE generation (as well as VRE capacities, see Figure S3 (a)) in the power sector, reduce VRE curtailment (see Figure S3 (b)) and reduce the need for dispatchable resources like CCGT and battery storage. This boosting effect grows as the share of electrolyzer in the H$_\text{2}$ supply chain increases, due to either increasing CO$_\text{2}$ price and/or increasing H$_\text{2}$ demand. In the deep decarbonization scenarios with total H$_\text{2}$ demands equivalent to 20\% and 80\% FCEV penetration, the total VRE generation increases by 12\% and 55\%, and the VRE curtailment reduces by 3\%-5\% (see Figure S3 (b)), respectively.

As a result of reduced need in dispatchable power resources and cheaper electrolytic H$_\text{2}$ production, we see cost savings from sector coupling in Fig. \ref{fig-synergy} (b), increasing with CO$_\text{2}$ price and H$_\text{2}$ demand and approaching 14\% of the total cost\footnote{The levelized costs of electricity and H$_\text{2}$ reduce by the same rate, if the cost savings are proportionally allocated between the power and H$_\text{2}$ sectors.} of the two sectors in the decoupled model in the deep decarbonization scenario. Sector-coupling also leads to greater CO$_\text{2}$ emissions reduction (up to 25\% of the emission of the case with \$0/tonne CO$_\text{2}$, \$300/kW electrolyzer, and 20\% FCEV penetration) than the case without coupling, owing to the increased penetration of VRE generation in the power sector. Both the cost savings and emission reductions from sector coupling increase with increasing H$_\text{2}$ demand (either from transportation or other end-uses such as heating or industrial sectors), more than doubling for the 96\% decarbonization scenario as H$_\text{2}$ demand quadraples. The CO$_\text{2}$ emission reductions are high in the cases with no carbon policy because of the large existing emissions in that case. As the CO$_\text{2}$ price increase from \$50/tonne to \$1000/tonne, the emission reductions benefit from sector coupling increase since individual decarbonization becomes more expensive within each sector.

\subsection{Storage or flexible demand}
G2P generators are not cost-competitive in the results discussed above, even in the deep decarbonization scenario, because of the relatively high capital cost and the additional efficiency losses incurred in supplying power rather than H$_\text{2}$. While they may become economically feasible in the future with economics of scale, technology innovations, and/or efficient deployment strategies (like sharing power conversion systems with electrolyzers or renewable plants), the results imply that the sectoral power exchange between the H$_\text{2}$ and power sectors could be highly imbalanced, as opposed to energy storage (charging and discharging typically of the same order of magnitude).

\begin{figure}[!htb]
	\centering
	\includegraphics[width=\textwidth]{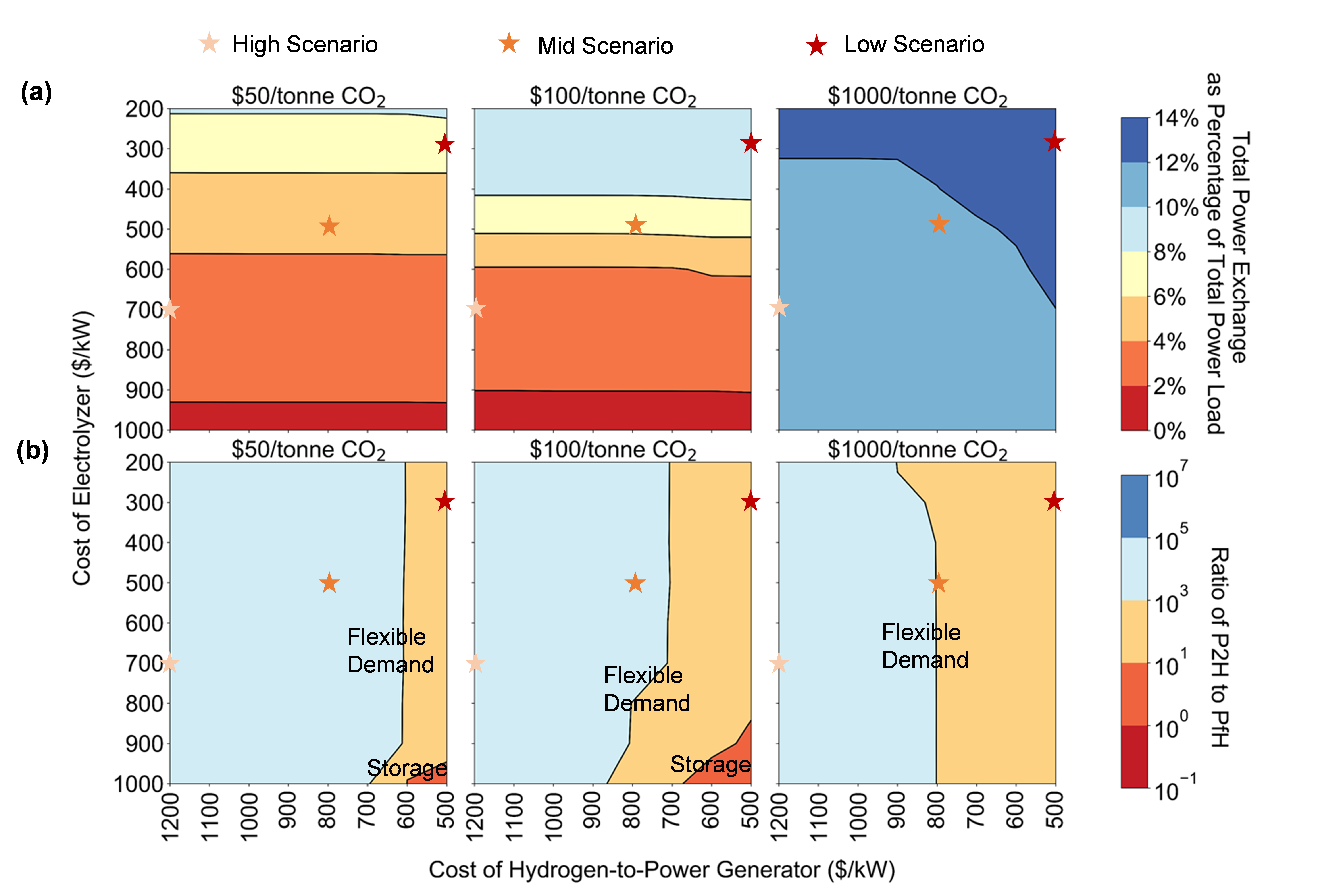}
	\caption{Power exchanges between H$_\text{2}$ supply chain and power system under different costs of electrolyzer and G2P generator. (a) Total power exchange throughput, power used to generate H$_\text{2}$ (P2H) plus power generated from H$_\text{2}$ (PfH), as percentages to the total power load (without H$_\text{2}$ generation); (b) The ratio of P2H to PfH. We define the role of H$_\text{2}$ supply chain as storage when the P2H and PfH are of the same order of magnitude, as flexible demand when the P2H is of higher order of magnitude than PfH, and as generator when the P2H is of lower order of magnitude than PfH. The cost scenario of electrolyzer ranges from \$1000/kW$_\text{e}$ to \$200/kW$_\text{e}$, and the cost scenario of H$_\text{2}$-based power generator ranges from \$1200/kW$_\text{e}$ to \$500/kW$_\text{e}$. The H$_\text{2}$ demand scenario corresponds to 20\% FCEV penetration. All other parameters are the same as the base case.}
	\label{fig-exchange}
\end{figure}

We further evaluate the role of H$_\text{2}$ supply chain and its interaction with the power sector for different values of electrolyzer and G2P generator capital costs in Fig. \ref{fig-exchange}. The electricity and H$_\text{2}$ interactions are quantified using two model outputs: a) the sum of annual power for H$_\text{2}$ production (P2H) and annual power generated from H$_\text{2}$ (PfH) exports as a percent of total annual electricity demand (without power for H$_\text{2}$ generation) and b) the ratio of annual P2H to PfH. The red, orange, and light pink stars in Fig. \ref{fig-exchange} represent high, medium and low capital costs for P2G and G2P technologies, respectively. We produce these scenarios based on potential cost savings from economies of scale, technology learning, and system designs that saves critical component costs like power electronics. We can see from Fig. \ref{fig-exchange} (a) that in the \$100/tonne CO$_\text{2}$ price case, along the cost-reduction pathway (i.e. a straight line connecting the three star markers), the power exchange between the H$_\text{2}$ and the power sectors increase from 3\% to 9\% of the total electric load. For the deep decarbonization scenario (\$1000/tonne CO$_\text{2}$), the power exchange can be as high as 14\% of the total electric load in the low technology cost scenario. Notably, as opposed to most literature that study H$_\text{2}$ as long-term/seasonal storage, our results in Fig. \ref{fig-exchange} (b) indicate that H$_\text{2}$ will more likely serve as a flexible demand response resource rather than long-term storage. This can be seen in Fig. \ref{fig-exchange} (b), where annual P2H generally tend to be several orders of magnitude greater than PfH, for almost the full technology cost space studied here. This observation stems from the additional efficiency losses and capital costs of converting H$_\text{2}$ back to power associated with H$_\text{2}$-based electricity storage vs. its use as a flexible demand resource. The only exceptions when the amounts of P2H and PfH are of the same order of magnitude under moderate decarbonization scenarios (\$50-100/tonne CO$_\text{2}$) with G2P capital costs below \$700/kW and electrolyzer capital costs greater than \$800/kW. That the minimum cost of G2P generator for the H$_\text{2}$ supply chain to play a storage role is higher in the moderate decarbonization scenario compared to the deep decarbonization scenario implies that the need for electrolytic H$_\text{2}$ increases faster than the need in PfH.

\subsection{Transmission deployment}
\begin{figure}[!htb]
	\centering
	\includegraphics[width=\textwidth]{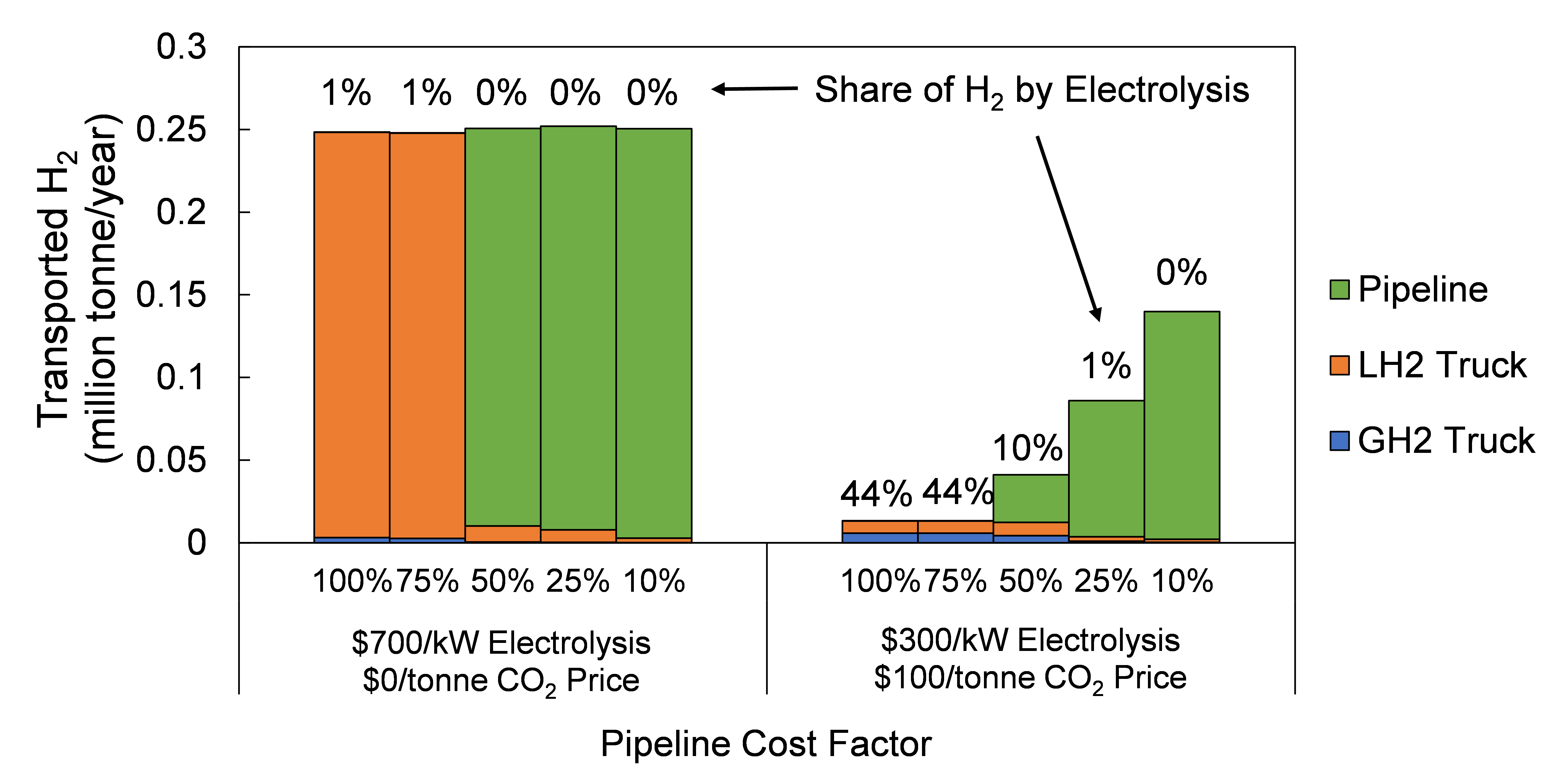}
	\caption{The amounts of transported H$_\text{2}$ per year via different transport modes under different pipeline cost scenarios and dominating H$_\text{2}$ generation modes. The FCEV penetration is 20\%. LH2: liquid H$_\text{2}$; GH2: Gaseous H$_\text{2}$.}
	\label{fig-transmission}
\end{figure}

In the results presented above, H$_\text{2}$ pipeline is not economically competitive compared to H$_\text{2}$ trucks. However, it could become attractive if the existing natural gas pipelines can be retrofitted at relatively low cost for transporting H$_\text{2}$, which has been a rising interest for the gas industry \cite{van2020hydrogen}. Therefore, we conduct a sensitivity analysis on the pipeline cost, given the uncertainty in the pipeline retrofitting cost as well as the cost of deploying dedicated H$_\text{2}$ pipelines. Besides pipeline capital cost, we focus on two scenarios, one with SMR as the dominating generation source, the other with significant electrolyzer generation while SMR playing as back up resource. As seen in Fig. \ref{fig-transmission}, H$_\text{2}$ pipeline becomes cost-competitive with H$_\text{2}$ truck if the pipeline capital cost is 50\% of the cost assumed in the base case. This can also be interpreted to represent a threshold cost for retrofitting existing natural gas pipelines to be compatible with 100\% H$_\text{2}$ flows. The other notable observation is that we need more H$_\text{2}$ transportation in the case with SMR as the dominant H$_\text{2}$ supply source compared to the case with significant electrolytic H$_\text{2}$ supply, and the share of electrolytic H$_\text{2}$ production goes down as pipelines become cheaper. This finding is, in part, driven by the assumption that distributed electrolyzers are allowed in the highly urbanized zone 4, while centralized SMRs are not.
Pipelines, because of their relatively high capital costs, are more cost-efficient to cope with large and steady H$_\text{2}$ transmission demand from centralized SMR production pathway, while electrolyzers can be deployed in a more distributed manner and thus complement smaller-scale and more flexible H$_\text{2}$ transmission mechanisms like trucks. The above conclusions also extended for increased H$_\text{2}$ demand scenarios, corresponding to FCEV penetration of 80\% (See Figure S4).

\section{Conclusions}
The interest in H$_\text{2}$ for decarbonizing energy systems is unquestionably increasing, in part driven by declining technology costs, greater policy emphasis on decarbonizing non-electric end-uses, and recognition of the limitations of direct electrification in certain applications. The unique versatility of H$_\text{2}$ as an energy carrier and its multiple uses, however, require a holistic view to accurately explore its role in future low-carbon energy systems and the accompanying technology pathways. Additionally, such a view could demonstrate the relative economic and environmental merits of H$_\text{2}$ and electricity use for various end-uses as well as their complementarity as vectors for decarbonizing the energy system. 

To this end, we developed a generalized framework for cost-optimal energy infrastructure investment and operations for decarbonizing multiple end-use sectors based on coordinated use of electricity and H$_\text{2}$ supply chains to manage spatio-temporal variations in renewable energy inputs and energy demands. This modelling approach provides numerous insights on the technological make-ups of these energy supply chains, spanning production, transport, storage and end-use, and their impacts on the cost of decarbonization, as highlighted via the U.S Northeast case study.  

First, in the coupled energy system, CCS is deployed at lower carbon prices in the H$_\text{2}$ sector than the power sector, which can be interpreted as CCS being more competitive in the H$_\text{2}$ supply chain than the power supply chain. This conclusion, however, goes counter to the observation that six times more CCS projects will be online before this decade in the power sector than for H$_\text{2}$ production\cite{Global_CCS_2020}. For regions like Europe, where decarbonization via H$_\text{2}$ is part of the government’s decarbonization roadmap, our study highlights the importance of prioritizing CCS deployment for H$_\text{2}$ production.

Second, power and H$_\text{2}$ sector coupling via flexible electrolysis and H$_\text{2}$ storage enables increased VRE penetration in the power sector, thereby reducing the need for alternative flexible resources for managing VRE variability (e.g. gas generation, battery storage, etc.) and in turn reducing total system cost.  Moreover, as opposed to other power-sector focused studies that emphasize H$_\text{2}$’s value as a grid-scale storage resource\cite{Guerra2020,Pellow15,Clerjon19}, our multi-sector view highlights the greater system value of P2G as a flexible demand resource that avoids the additional efficiency losses and capital cost incurred with P2G2P pathways. This conclusion is found to be robust to future expectations on the capital costs of electrolyzer and G2P systems. 
Since electrolyzers and H$_\text{2}$ storage are commercially available, this finding also suggests that H$_\text{2}$ playing a role for grid balancing could be sooner than the full P2G2P routes becoming cost-effective. 

Third, as compared to the independent optimization of each supply chain, we find that sector-coupling via P2G, reduces the cost of energy system decarbonization and that this benefit grows as the demand for H$_\text{2}$ in other end-use sectors increases. Realizing the benefits of such cross-sector coordination, however, calls for policy and market reforms. For example, H$_\text{2}$ prices need to be settled at similar spatiotemporal resolution as electricity prices, to provide incentives and signals for H$_\text{2}$ infrastructure owners, and electrolyzers should be allowed to provide ancillary services to power systems. Moreover, both integrated operation and planning of power and H$_\text{2}$ sectors, through a shared independent system operator, could help fully exploit the sector-coupling benefits.

Finally, we found that the choice of H$_\text{2}$ transport infrastructure is intricately dependent on the choice of H$_\text{2}$ production infrastructure, with pipelines more synergistic with natural gas based production pathways because of the matching scale of production and transmission capacity of the two assets.

There are a number of areas for future work that can build on this analysis. While this study has focused on H$_\text{2}$ used in transportation, greater H$_\text{2}$ demand might be realized from decarbonizing heating and industrial sectors, such as ammonia production and steel manufacturing. Those H$_\text{2}$ demands will have different temporal profiles and flexibility compared to FCEV charging, and thus may affect the H$_\text{2}$ supply-demand balance in different ways. Accounting for these heterogeneous H$_\text{2}$ demands should be further explored. Investigating the last-mile delivery of H$_\text{2}$ is out of scope of this study but is definitely a key area for future work. The last-mile H$_\text{2}$ distribution network will affect the total number of gas trucks needed for H$_\text{2}$ transport in the system. Increased traffic congestion and safety considerations of a very large H$_\text{2}$ truck fleet may limit the deployment of trucks and increase the values of pipeline based transport, and further analysis in conjunction with traffic simulations is needed to understand this aspect. 

Regional factors, including resource availability and demand level, could significantly affect the optimal technology portfolios and the costs in both power and H$_\text{2}$ sectors. For example, while natural gas supply is abundant in the studied U.S. Northeast region, it may be insufficient or expensive for many regions in the world, leading to higher shares of VRE, electrolytic H$_\text{2}$ production, and storage. The ability to cost-effectively store H$_\text{2}$ at various time scales and capacities is a critical factor in determining the optimal system architecture. Underground caverns, where available, can provide cheap H$_\text{2}$ storage for VRE and electrolysis deployment.  The competing role of CCS with H$_\text{2}$ storage in underground resources also need to be considered for relevant regions.

\section*{Conflicts of interest}
There are no conflicts to declare.

\section*{Acknowledgements}
This work was partially supported by Shell New Energies Research and Technology, Amsterdam, Netherlands, and the Low-Carbon Energy Centers on Electric Power Systems and Carbon Capture Sequestration and Utilization at MIT Energy Initiative. We would like to thank Dr. Joe Powell, Dr. Mark Klokkenburg, and Prof. Robert Armstrong for their valuable advice on this work. We acknowledge the MIT SuperCloud and Lincoln Laboratory Supercomputing Center for providing resources that have contributed to the research results in this paper.


\bibliography{Reference}

\end{document}


\title{Sector coupling via hydrogen to lower the cost of energy system decarbonization: Supplementary Information}
\author{Guannan He\footnotemark[1]~~\footnotemark[2]~~,
    Dharik S. Mallapragada\footnotemark[2]~~,
	Abhishek Bose\footnotemark[2]~~,\\
	Clara F. Heuberger\footnotemark[3]~~, 
	Emre Gen\c{c}er\footnotemark[2]
}
\date{}
\thispagestyle{empty}
\pagestyle{empty}
\maketitle

\footnotetext[1]{Email: gnhe@mit.edu}
\footnotetext[2]{MIT Energy Initiative, Massachusetts Institute of Technology, Cambridge, MA, USA.}
\footnotetext[3]{Shell Global Solutions International B.V., Shell Technology Centre Amsterdam, 1031 HW Amsterdam, Netherlands.}

\section*{Supplementary Figures}

\begin{figure}[H]
    \centering
    \includegraphics[width = \textwidth*2/3]{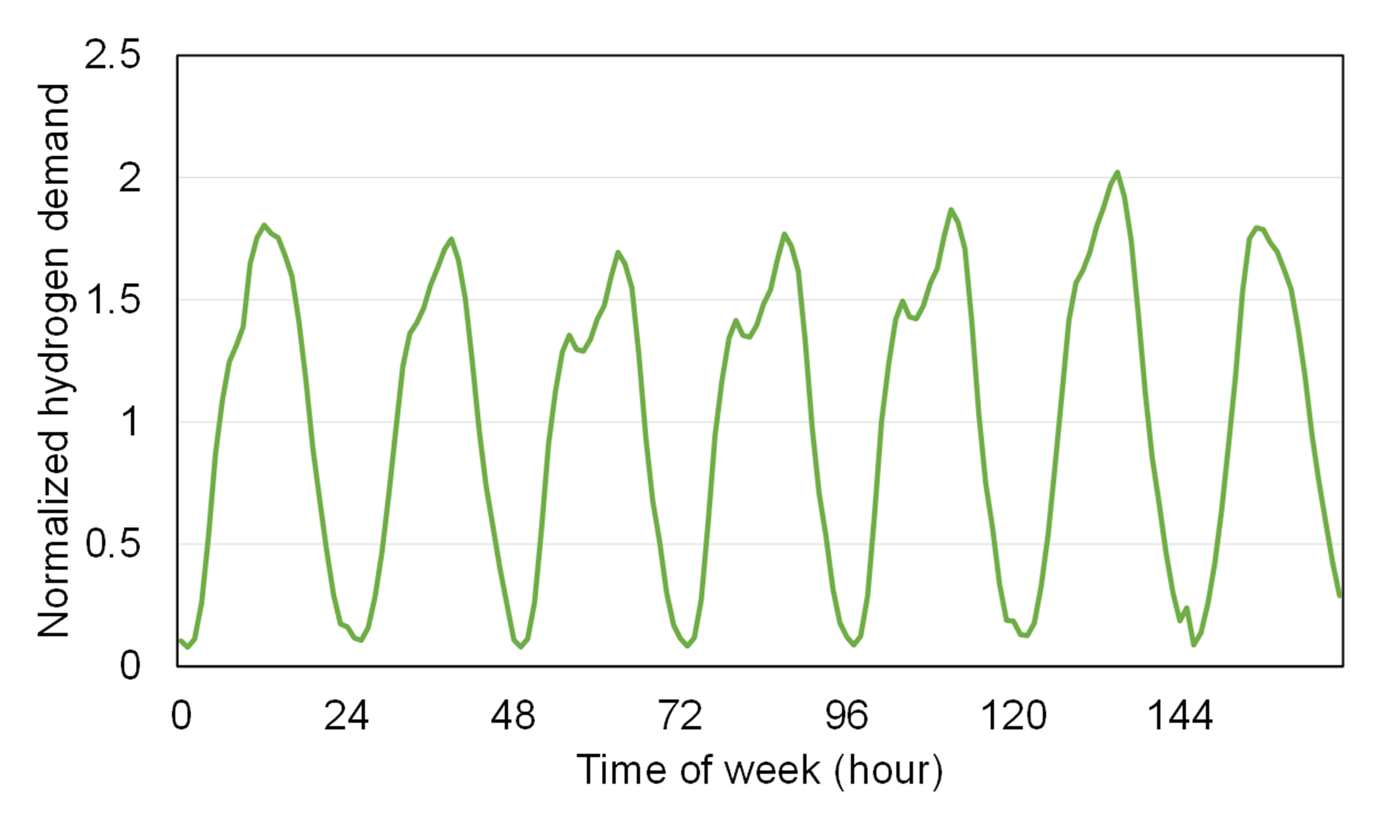}
    \caption{Hourly hydrogen refuelling profile normalized based on the mean}
    \label{fig:Unit Refueling Profile}
\end{figure}

\begin{figure}[H]
    \centering
    \includegraphics[width = \textwidth]{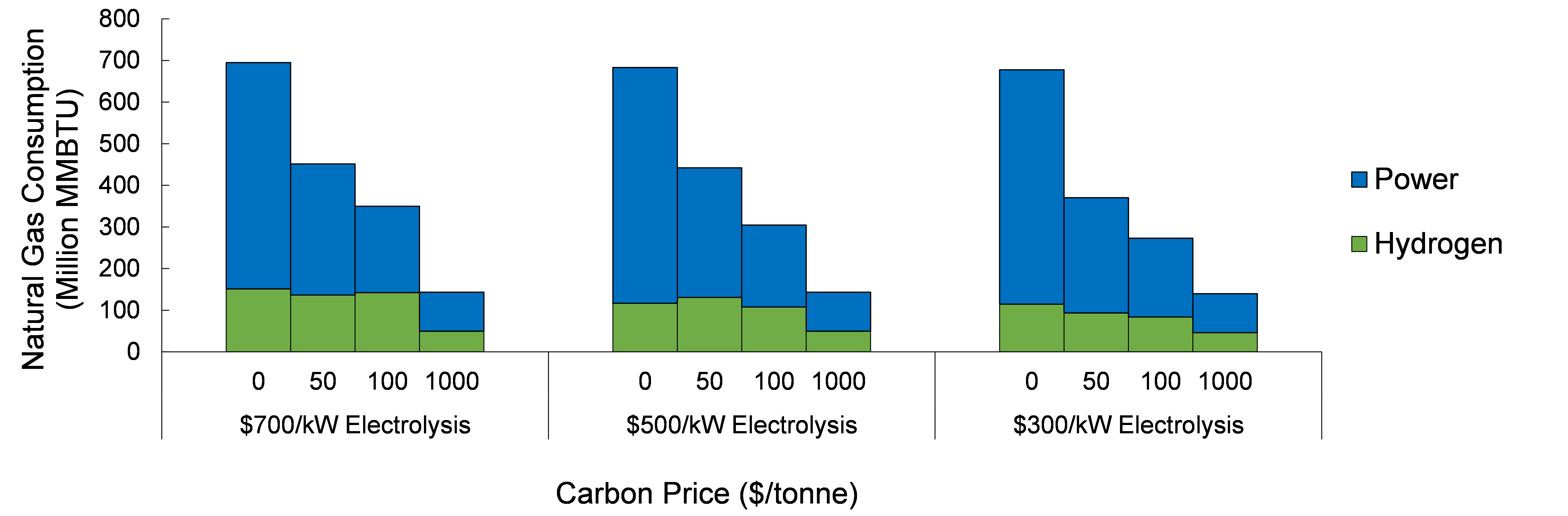}
    \caption{Natural gas consumption breakdowns in H$_2$ and power sectors.}
    \label{fig:NG Consumption}
\end{figure}

\begin{figure}[H]
    \centering
    \includegraphics[width = \textwidth]{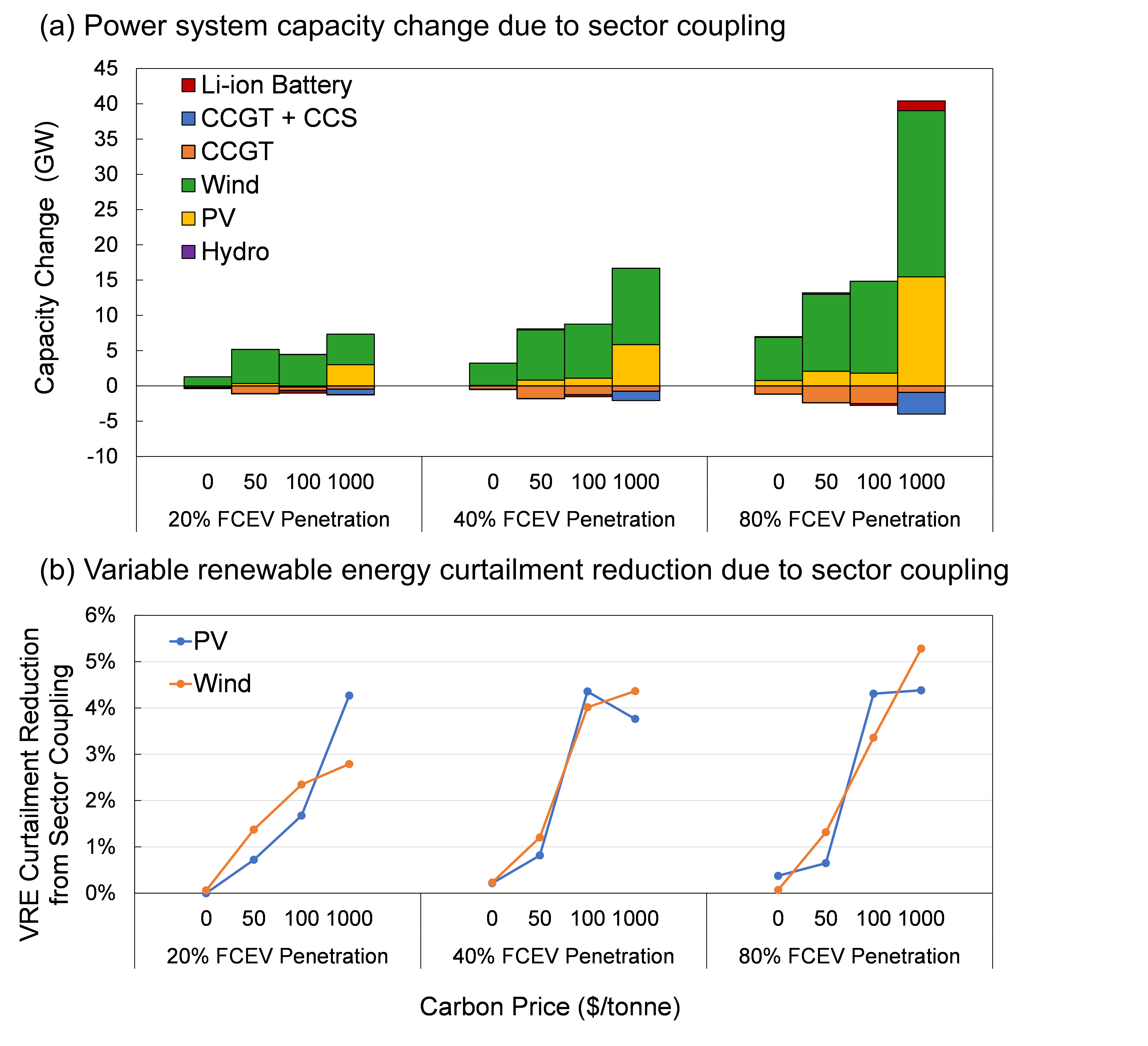}
    \caption{Differences in the optimal power system capacity mix and variable renewable energy (VRE) curtailment between energy systems with and without conversions between power to H$_\text{2}$ under various CO$_\text{2}$ price and FCEV penetration scenarios. (a) Differences in the optimal power system capacity mix; (b) Differences in the VRE curtailment. In (b), the VRE curtailment reductions are calculated as the differences in the curtailment percentages between energy systems with and without conversions between power to H$_\text{2}$ in each CO$_\text{2}$ price and FCEV penetration scenario.}
    \label{fig:synergy_power_cap}
\end{figure}

\begin{figure}[H]
    \centering
    \includegraphics[width = \textwidth]{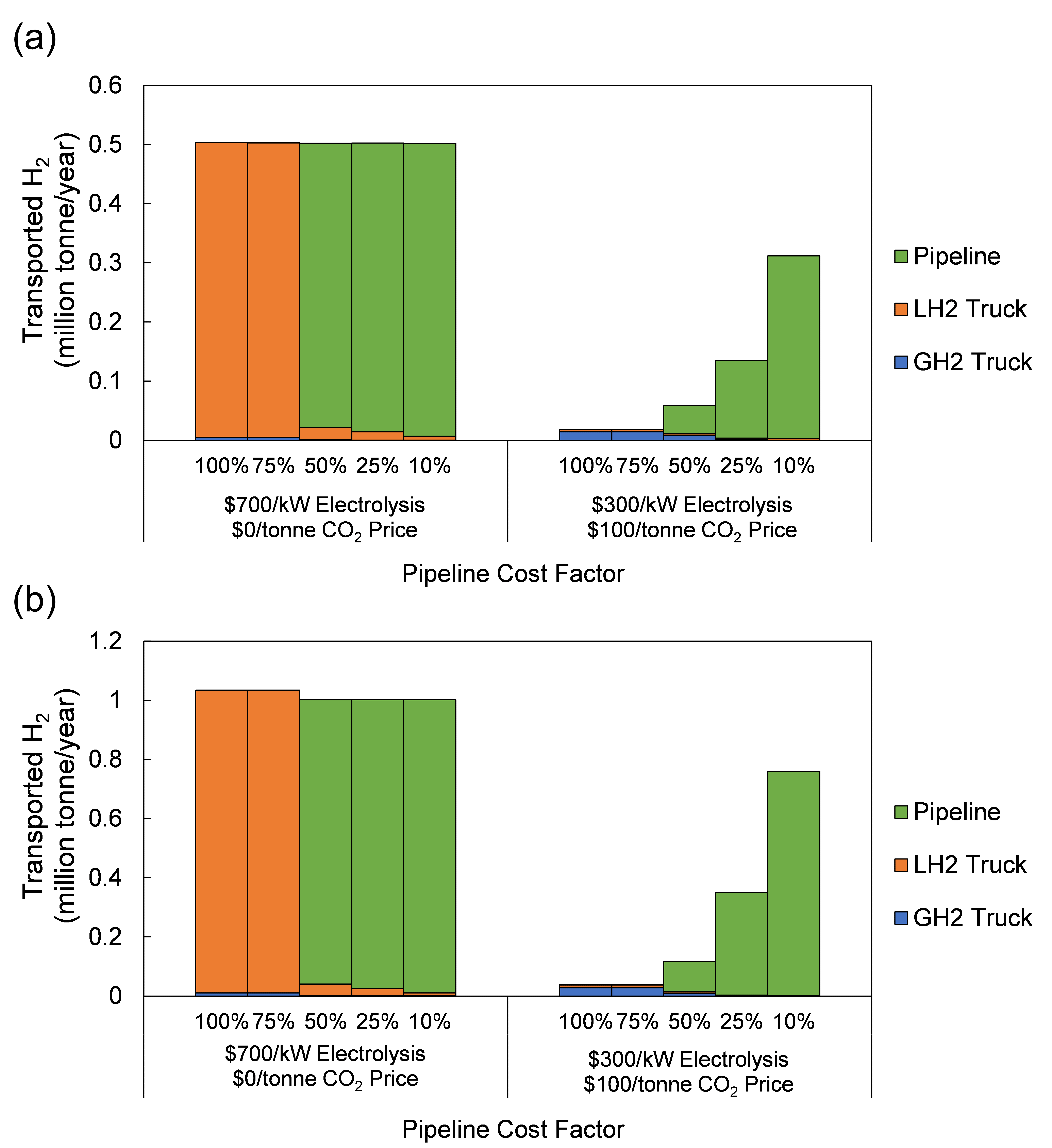}
    \caption{The amounts of transported H$_2$ per year via different transport modes under different pipeline cost scenarios and dominating H$_2$ generation modes. (a) 40\% FCEV penetration; (b) 80\% FCEV penetration. LH2: Liquid H$_2$; GH2: Gaseous H$_2$.}
    \label{fig:transmission_SI}
\end{figure}

\newpage
\section*{Supplementary Tables}

\begin{table}[h]
  \centering
  \caption{Additional Parameters of the e-H$_2$ model}
    \begin{tabular}{lc}
    \toprule
    Discount rate & 7\% \\
    Power transmission expansion cost & 1,600/MW-mile \\
    Power transmission loss & 1\%/100 miles \\
    Value of lost load (electricity) & \$20,000/MWh \\
    Value of lost load (hydrogen) & \$1,000/kg \\
    Gas Price  & \$5.4/MMBTU \\
    CO$_2$ transportation cost & \$20/tonne \\
    \bottomrule
    \end{tabular}%
  \label{tab:add_par}%
\end{table}%

\begin{table}[htbp]
  \centering
  \caption{Parameters of Existing Power Transmission Lines. Zone 7 represents Canada.}
    \begin{tabular}{lcc}
    \toprule
    Line  & Max Power Flow (MW) & Line Distance (km) \\
    \midrule
    Zone 1-2 & 2000  & 317 \\
    Zone 2-3 & 2950  & 199 \\
    Zone 3-4 & 760   & 99 \\
    Zone 4-5 & 1528  & 216 \\
    Zone 3-5 & 600   & 158 \\
    Zone 2-5 & 800   & 179 \\
    Zone 5-6 & 5400  & 186 \\
    Zone 2-6 & 150   & 340 \\
    Zone 6-7 & 2600  & 0 \\
    Zone 2-7 & 1650  & 0 \\
    Zone 1-7 & 800   & 0 \\
    \bottomrule
    \end{tabular}%
  \label{tab:addlabel}%
\end{table}%

\begin{table}[H]
  \centering
  \caption{Inter-Zone Distances of the Six Zones in the U.S. Northeast for Trucks and Pipelines (mile).}
  \vspace{1em}
    \begin{tabular}{c|cccccc}
    \toprule
    \textbf{Zone} & \textbf{1} & \textbf{2} & \textbf{3} & \textbf{4} & \textbf{5} & \textbf{6} \\
    \hline
    \textbf{1} & 0   & 317   & 504   & 602   & 487   & 608 \\
    \textbf{2} & 317   & 0     & 199   & 297   & 179   & 340 \\
    \textbf{3} & 504   & 199   & 0     & 99    & 158   & 333 \\
    \textbf{4} & 602   & 297   & 99    & 0     & 216   & 358 \\
    \textbf{5} & 487   & 179   & 158   & 216   & 0     & 186 \\
    \textbf{6} & 608   & 340   & 333   & 358   & 186   & 0 \\
    \bottomrule
    \end{tabular}%
  \label{tab:zone_dist}%
\end{table}%

\newpage




\printbibliography[title=Supplementary Reference]